\documentclass[12pt,reqno]{amsart} 

\usepackage{mathrsfs,amsfonts,amsthm,amsmath,amssymb,thmtools} 
\usepackage[dvipsnames]{xcolor} 
\usepackage{extarrows}

\usepackage{enumitem}                
\usepackage{lineno}                 

\usepackage[numbers,sort&compress]{natbib}  

\usepackage{graphics}
\usepackage{graphicx}   
\graphicspath{{images/}} 
\usepackage{subfig}     
\usepackage{rotating}

\usepackage[hyperindex, pdfencoding = auto, psdextra, bookmarksdepth = 4]{hyperref}  
\hypersetup{
    bookmarksopen      = true,
    bookmarksopenlevel = 1,
    bookmarksnumbered  = true,
    pdftitle    = {},
    pdfauthor   = {},
    linktoc     = page,
    anchorcolor = green,
    colorlinks  = true,        
    linkcolor   = black,
    citecolor   = blue,
    filecolor   = blue,
    urlcolor    = magenta,
}

\newcommand{\NN}{\mathbb{N}}

\newcommand{\RR}{\mathbb{R}}

\newcommand{\cA}{\mathcal{A}}
\newcommand{\cB}{\mathcal{B}}

\newcommand{\cD}{\mathcal{D}}
\newcommand{\cE}{\mathcal{E}}
\newcommand{\cF}{\mathcal{F}}

\newcommand{\cH}{\mathcal{H}}
\newcommand{\cI}{\mathcal{I}}

\newcommand{\cM}{\mathcal{M}}
\newcommand{\cN}{\mathcal{N}}

\newcommand{\cP}{\mathcal{P}}

\newcommand{\cR}{\mathcal{R}}

\usepackage[capitalise, noabbrev, nameinlink]{cleveref}    

\theoremstyle{plain}

\newtheorem{theorem}{Theorem}[section]
\newtheorem{lemma}[theorem]{Lemma}

\newtheorem{proposition}[theorem]{Proposition}
\newtheorem{corollary}[theorem]{Corollary}

\theoremstyle{definition}   

\newtheorem{remark}[theorem]{Remark}

\makeatletter
    \def\subsection{\@startsection{subsection}{2}%
    \z@{.5\linespacing\@plus.7\linespacing}{.3\linespacing}%
    {\normalfont\bfseries}}
\makeatother 

\makeatletter                                         
    \newcommand{\LeftEqNo}{\let\veqno\@@leqno}        
\makeatother

\makeatother 

\numberwithin{equation}{section}                      

\begin{document}

\

\vspace{-2cm}

\title[Orthogonal lifting and amenability in ultrapowers]
{Orthogonal lifting and amenability in ultrapowers}

\author{Minghui Ma}
\address{Minghui Ma, School of Mathematical Sciences, Dalian University of Technology, Dalian, 116024, China}
\email{minghuima@dlut.edu.cn}

\author{Wenming Wu}
\address{Wenming Wu, School of Mathematical Sciences, Chongqing Normal University, Chongqing, 401331, China}
\email{wuwm@amss.ac.cn}

\author{Yifeng Ye}
\address{Yifeng Ye, School of Mathematical Sciences, Chongqing Normal University, Chongqing, 401331, China}
\email{yyf@amss.ac.cn}

\keywords{Orthogonal
lifting, amenability, Connes embeddability,
type $\mathrm{II}_1$ factor,
ultrapower}

\makeatletter
\@namedef{subjclassname@2020}{\textup{2020} Mathematics Subject Classification}
\makeatother
\subjclass[2020]{Primary 46L10; Secondary 46M07}

\begin{abstract}
	Let $\cM$ be a  type $\mathrm{II}_1$
	factor  and let 
	 $\omega$ be  a free ultrafilter on $\NN$.
	We prove that any
	two  orthogonal separable
	amenable von Neumann subalgebras in  the ultrapower
	$\cM^{\omega}$ admit orthogonal finite-dimensional  
	lifts. In particular, two
	independent unitary operators
	in $\cM^\omega$ admit independent
	finite-spectrum representatives.
	We also show that 
	the orthogonal  lifting can fail  whenever  one of
	the subalgebras is  nonamenable and Connes-embeddable.
	\end{abstract}
	\maketitle


\section{Introduction}
 The seminal papers of McDuff \cite{McD70} and Connes \cite{Co76}
in the 1970s established ultraproducts and ultrapowers of
tracial von Neumann algebras as important tools in the study
of von Neumann algebras.
Since ultraproducts and ultrapowers are constructed as limits
along a free ultrafilter, a natural and challenging problem
is to choose representing sequences that preserve specified
algebraic relations or independence properties exactly
at each component; see \cite{CIKE23,HI24,PS26}.
The author in \cite[Problem~30]{KE25}
 asks whether any pair of independent unitary  operators 
in the ultrapower of a separable 
type $\mathrm{II}_1$ factor can be lifted to sequences of 
unitary operators that are independent at each component.

 Let $(\cM,\tau)$ be a tracial von Neumann algebra and
let $\cA,\cB\subseteq\cM$ be von Neumann subalgebras.
We say that $\cA$ and $\cB$ are \emph{$n$-independent} if
\[
\tau(z_1\cdots z_k)=0
\] 
whenever $1\leqslant k\leqslant 2n$ and 
$z_1,\ldots,z_k$ are trace-zero operators taken alternately
from $\cA$ and $\cB$.
 Popa
\cite[Section~3]{Po14KS} used finite-order independence
to study  the  Kadison--Singer problem.
If $\cA$ and $\cB$ are $n$-independent for all $n\geqslant 1$, we say that $\cA$ and $\cB$ are \emph{freely
independent} or \emph{free} in short.
 This agrees with Voiculescu's notion of free independence
\cite{Vo85}.
 Two von Neumann subalgebras $\cA,\cB\subseteq\cM$
are called \emph{independent}, or \emph{orthogonal},
written as $\cA\perp_\tau\cB$, if they are $1$-independent;
equivalently,
\[
 \tau(ab)=\tau(a)\tau(b)
\qquad(a\in\cA,\ b\in\cB).
\] 
Two unitary operators $u,v\in \cM$ are \emph{independent} if
$W^*(u)\perp_\tau W^*(v)$.

Chifan, Ioana, and Kunnawalkam Elayavalli
\cite[Lemma~3.1]{CIKE23} proved that orthogonal abelian
von Neumann subalgebras of dimensions two and three
in an ultraproduct of type $\mathrm{II}_1$ factors admit
orthogonal lifts.
This result was a key ingredient in their construction
of a non-Gamma type $\mathrm{II}_1$ factor not elementarily
equivalent to any free group factor.
Houdayer and Ioana
\cite[Theorem~5.1 and the subsequent discussion]{HI24}
subsequently showed that any two separable abelian
$2$-independent von Neumann subalgebras of a tracial
ultraproduct admit orthogonal abelian lifts.
They combined this lifting theorem with asymptotic
freeness to construct a non-Gamma type $\mathrm{II}_1$ factor
not elementarily equivalent to any free product of
diffuse tracial von Neumann algebras.
More recently, Patchell and Shiner
\cite[Theorem~C]{PS26} proved that any two separable
abelian $2n$-independent von Neumann subalgebras of an
ultraproduct admit $n$-independent
finite-dimensional abelian lifts.
These results provide   affirmative answers to special
cases of Problem~30 in \cite{KE25}.

Throughout this paper, $\omega$ denotes a free ultrafilter on $\NN$.
A von Neumann algebra is called \emph{separable} if it has a
separable predual.
Recall that a separable tracial von Neumann algebra $\cN$
is \emph{amenable} if and only if it contains an increasing
sequence of finite-dimensional von Neumann subalgebras
whose union is weakly dense in $\cN$; see \cite{Co76}.
Our first theorem shows that any two orthogonal separable 
  amenable von Neumann subalgebras 
  of  $\mathcal{M}^\omega$ 
   admit orthogonal
   finite-dimensional lifts. 

\begin{theorem}\label{thm:amenable}
 Suppose that $(\cM,\tau)$ is a type $\mathrm{II}_1$ factor,
 $\cA,\cB$ are   separable
amenable von Neumann
subalgebras of  the ultrapower $\cM^{\omega}$
and $\cA\perp_{\tau_\omega} \cB$. Then there are
finite-dimensional von Neumann subalgebras $\cA_n,\cB_n\subseteq \cM$ satisfying
\begin{equation*}
 \cA_n \perp_{\tau} \cB_n\quad(n\geqslant 1),\qquad
 \cA \subseteq\prod_{\omega} \cA_n,\qquad
 \cB \subseteq\prod_{\omega} \cB_n.
\end{equation*}
Moreover, if $\cA$ $($respectively, $\cB)$ is abelian, then the algebras
$\cA_n$ $($respectively, $\cB_n)$
can be chosen to be abelian.
\end{theorem}

As an application of \Cref{thm:amenable}, we  obtain an 
affirmative answer to \cite[Problem~30]{KE25}; see \Cref{cor:unitaries}.

Recall that a  separable tracial von Neumann algebra
$(\cN,\tau_{\cN})$ is called \emph{Connes-embeddable}  if it admits a
unital trace-preserving $^*$-embedding into $\cR^\omega$,
where $\cR$ denotes the hyperfinite type $\mathrm{II}_1$ factor;
see \cite{Co76}.
Jung \cite{Ju07} proved that a finitely generated
Connes-embeddable  tracial von Neumann algebra is amenable
if and only if any two trace-preserving $^*$-embeddings into
$\cR^\omega$ 
 are conjugate by a unitary operator in $\mathcal{R}^\omega$.
Atkinson and Kunnawalkam Elayavalli
\cite[Theorem~4.7]{AKE21} extended Jung's characterization
to separable Connes-embeddable tracial von Neumann algebras 
by
showing that such an algebra $\cN$ is amenable if and only if
any two trace-preserving $^*$-embeddings 
$\cN\to\cR^\omega$
are conjugate via a sequence of trace-preserving unital
completely positive maps on $\cR$.
This
characterization will be used in 
 the following theorem,
 which shows that the
 amenability assumption in \Cref{thm:amenable} is essential.

\begin{theorem}\label{thm:counterexample}
  Let $(\cM,\tau)$ be a type $\mathrm{II}_1$ factor, and let
 $(\cN,\tau_{\cN})$ be a separable nonamenable
  Connes-embeddable tracial von Neumann algebra.
  For every nontrivial separable  von Neumann subalgebra
  $\cB\subseteq\cM^\omega$ and every nontrivial 
  projection
  $q\in\cB$, there exists a  von Neumann
  subalgebra $\cA\subseteq\cM^\omega$ satisfying

\begin{enumerate}[label=(\arabic*), leftmargin=*]
  \item[$(1)$]   $(\cA,\tau_\omega|_{\cA})\cong(\cN,\tau_{\cN})$
  as tracial von Neumann algebras;
  \item[$(2)$]  $\cA\perp_{\tau_\omega}\cB$;
  \item[$(3)$] $\cA$ and $W^*(q)$ admit no orthogonal
  lifts; that is, there  exist no sequences of   
  von Neumann subalgebras
  $\cA_n,\cD_n\subseteq\cM$ such that
  \[
    \cA_n \perp_\tau \cD_n\quad(n\geqslant 1),\qquad
  \cA\subseteq\prod_\omega\cA_n,\qquad
  W^*(q)\subseteq\prod_\omega\cD_n.
  \]
\end{enumerate} 
  \end{theorem}

  Fix a nontrivial separable amenable von Neumann subalgebra
$\cB\subseteq\cM^\omega$.
By \Cref{thm:amenable} and \Cref{thm:counterexample},
a separable Connes-embeddable tracial von Neumann algebra
$(\cN,\tau_{\cN})$ is amenable if and only if
$\varphi(\cN)$ and $\cB$ admit orthogonal lifts for every
trace-preserving $^*$-embedding $\varphi:\cN\to\cM^\omega$
satisfying $\varphi(\cN)\perp_{\tau_\omega}\cB$.
 
This paper is organized as follows. 
In \Cref{sec:prelim}, we recall the construction of ultraproducts
and prove some preliminary lemmas.
In \Cref{sec:amen}, we first construct orthogonal lifts for pairs of
finite-dimensional subalgebras and then use hyperfinite approximation
to prove \Cref{thm:amenable}. As a consequence, we show that independent
unitary operators  have finite-spectrum representatives that are independent at
every component. 
In \Cref{sec:counterexamples}, we construct orthogonal
subalgebras of $\cM^\omega$ that admit no orthogonal lifts.

\section{Preliminaries}\label{sec:prelim}
Let $\{ (\cM_n, \tau_n)\}_{n \geqslant 1}$ be a sequence of tracial von Neumann algebras.  
The product algebra $\prod_{n \geqslant 1} \cM_n$ 
is the set of bounded sequences 
$x=\left(x_n\right)$ with  
 $x_n \in \cM_n$ for 
each $n\geqslant 1$, endowed with the norm $\|x\|=\sup_n \left\|x_n\right\|$.
 The  (\emph{tracial}) \emph{ultraproduct} $\prod_\omega \cM_n$ 
 is the quotient of $\prod_{n=1}^\infty \cM_n$ by 
 the norm-closed two-sided ideal 
 $$\cI_\omega:=\left\{ (x_n) \in \prod_{n \geqslant 1} \cM_n:\lim_{n \to \omega} \tau_n\left(x_n^* x_n\right)=0 \right\}.$$
 For $(x_n) \in\prod_{n\geqslant1}\cM_n$, we denote its
equivalence class in $\prod_\omega\cM_n$ by $(x_n)_\omega$.
Equipped with the trace
 $$\tau_\omega ((x_n)_\omega)=\lim_{n \to \omega} \tau_n(x_n)$$
on $\prod_\omega \cM_n$, $(\prod_\omega \cM_n,\tau_\omega)$ 
is a tracial von Neumann algebra.

When these  $(\cM_n, \tau_n)$ are the same tracial von Neumann algebra $(\cM, \tau)$, 
we set $\cM^\omega=\prod_\omega \cM$, and we say that  $(\cM^\omega, \tau_\omega)$  
is the  (\emph{tracial}) \emph{ultrapower} of $(\cM, \tau)$ along $\omega$.
Moreover, 
if  $\cM$  is a type $\mathrm{II}_1$ factor,  
then $\cM^\omega$ is also a  
type $\mathrm{II}_1$ factor, see  \cite[Lemma~A.4.2]{SS08}.
The reader is referred to \cite{GH01,AH14,FGL06} for 
a general study of ultraproducts of
operator algebras.

Let $\cM$ be a type $\mathrm{II}_1$ factor. For any two unital $^*$-embeddings \mbox{$\pi_1,\pi_2: M_n(\mathbb{C})\rightarrow \cM$}, there is a unitary $u\in\cM$ such that $\pi_2(x)=u^*\pi_1(x)u$ for all $x\in M_n(\mathbb{C})$.
 Hence for any finite-dimensional von Neumann subalgebra $\cA$ of 
 $\cM^{\omega}$, 
 there is a unitary operator $u=(u_n)_\omega \in \cM^{\omega}$ such that $u^*\cA u\subseteq \cM\subseteq\cM^{\omega}$ with respect to the diagonal embedding of $\cM$.
 Conjugating the fixed copy of $\mathcal A$ in $\mathcal M$ by $u_n$ at each component yields the following lifting lemma.

\begin{lemma}\label{lem:fd-lifts}
  Let $(\cM,\tau)$ be a type $\mathrm{II}_1$ factor,
  and let 
  $\cF$ be a 
     finite-dimensional von Neumann subalgebra   of   
  $\cM^{\omega}$. Then there exist
    trace-preserving $^*$-embeddings
	\mbox{$\theta_n\colon \cF\to\cM$} such that $x=(\theta_n(x))_{\omega}$
  for each $x\in \cF$.
  \end{lemma}

To lift two independent separable amenable von Neumann subalgebras in $\cM^{\omega}$, we need the following 
closedness lemma.  
For a von Neumann algebra $\cM$, let $(\cM)_1$ denote its closed unit ball in the operator norm.

\begin{lemma}\label{lem:closed-lifts}
Let $(\cM,\tau)$ be a tracial von Neumann algebra and $m\geqslant1$.
Suppose that $S_n\subseteq\{(x_1,\ldots,x_m):  x_j \in (\cM)_1, \,1 \leqslant j \leqslant m\}$
 is nonempty for every $n\geqslant 1$.
Then the set
\begin{equation*}
 \prod_{\omega}S_n
 =\{(x_n)_{\omega}\colon x_n\in S_n\text{ for every }n\geqslant 1\}
 \subseteq(\cM^{\omega})^m
\end{equation*}
is closed in componentwise $2$-norm.
\end{lemma}

\begin{proof}

	For $x=(x_1,\ldots,x_m)\in\cM^m$, define
	\[
	d_n(x)
	:=\inf_{(y_1,\ldots,y_m)\in S_n}
	   \max_{1\leqslant j\leqslant m}\|x_j-y_j\|_2.
	\]
	The triangle inequality gives
	\begin{align}\label{sjbdsdn}
	|d_n(\xi)-d_n(\eta)|
	\leqslant
	\max_{1\leqslant j\leqslant m}\|\xi_j-\eta_j\|_2
	\qquad(\xi,\eta\in\cM^m).
\end{align}
	For $x=(x_1,\ldots,x_m)\in(\cM^\omega)^m$,
write $x_j$ as  $x_j=(x_{j,n})_\omega$
	and define
	\[
	d(x):=\lim_{n\to\omega}
		   d_n(x_{1,n},\ldots,x_{m,n}).
	\]
It follows from   \eqref{sjbdsdn}  
 that $d$ is independent
	of the choice of the representatives of $x$. 
	Moreover, for
	$x,y\in(\cM^\omega)^m$, \eqref{sjbdsdn}   yields
	\[
	|d(x)-d(y)|
	\leqslant
	\lim_{n\to\omega}
	\max_{1\leqslant j\leqslant m}\|x_{j,n}-y_{j,n}\|_2
	=
	\max_{1\leqslant j\leqslant m}\|x_j-y_j\|_2.
	\]
	Thus $d$ is continuous in the topology induced by
	the maximum of the component $2$-norms.
	
	Now we show that  $d^{-1}(\{0\})=\prod_\omega S_n$.
  Suppose that $d(x)=0$. For each $n \geqslant 1$,
	choose $(y_{1,n},\ldots,y_{m,n})\in S_n$ such that
	\[
	\max_{1\leqslant j\leqslant m}\|x_{j,n}-y_{j,n}\|_2
	\leqslant
d_n(x_{1,n},\ldots,x_{m,n})+\frac1n.
	\]
	Since $S_n\subseteq((\cM)_1)^m$, each sequence
$(y_{j,n})_n$ is uniformly bounded. Taking ultralimits
in the preceding inequality, we 
obtain   $x_j=(y_{j,n})_\omega$ for every
$1\leqslant j\leqslant m$, and hence
$x\in\prod_\omega S_n$.
Consequently, $\prod_\omega S_n=d^{-1}(\{0\})$
is closed in the componentwise $2$-norm topology.
	\end{proof}

The following lemma is an application of 
Popa's Theorem \cite[Theorem~0.1(b)]{Po14}.

\begin{lemma}\label{thm:popa}
Let $\cM$ be a type $\mathrm{II}_1$ factor, 
and let $\cA\subseteq\cM^{\omega}$ be a norm-separable
unital $C^*$-subalgebra. Then there is a  trace-zero 
symmetry operator
$v\in\cM^{\omega}$ such that   $\cA$ and $W^*(v)$ are free.
\end{lemma}
\begin{proof} As $\cA$ is norm-separable and $\cM^{\omega}$ 
 is nonseparable, there is a projection $p\in\cM^{\omega}\setminus\cA$ with $\tau_{\omega}(p)=\frac{1}{2}$. By \cite[Theorem~0.1(b)]{Po14}, there is a unitary 
  $u\in \cM^{\omega}$ such that $\cA''\vee u^*W^*(p)u \cong\cA''* W^*(u^*pu)$.
 Let $v=2u^*pu-1$. Then $v$ is the required symmetry operator.
\end{proof}

\section{Orthogonal lifting for amenable subalgebras}\label{sec:amen}

 In this section, we  prove \Cref{thm:amenable}. We begin with the notation needed
for the finite-dimensional case.

Let  $(\cM,\tau)$ be a  
type $\mathrm{II}_1$ factor,  and 
let $\cA_0$ and $\cB_0$  be  
nontrivial finite-dimensional von Neumann subalgebras 
of $\cM^\omega$ 
with $\cA_0 \perp_{\tau_\omega} \cB_0$.
Suppose that $p=(p_i)_{i=1}^r$ and \mbox{$q=(q_j)_{j=1}^s$}
are tuples of projections such that
$\{1,p_1,\ldots,p_r\}$ and $\{1,q_1,\ldots,q_s\}$
form linear bases of $\cA_0$ and $\cB_0$ respectively.
Write $C^*(p,q)$ for the unital $C^*$-algebra 
 generated by these projections.
Define a real linear map
$h_{p,q}\colon\RR^{r\times s}\to\cM^\omega$ by
\begin{equation}\label{eq:operator}
h_{p,q}(\xi)=\mathbf{i}\sum_{i=1}^r\sum_{j=1}^s\xi_{ij}[p_i,q_j]     \qquad  \big( \xi=(\xi_{ij})\in \RR^{r\times s} \big).
\end{equation}
Here $\mathbf{i}=\sqrt{-1}$. Define maps
$p_i(\cdot)\colon\RR^{r\times s}\to\cM^\omega$ and
$\Phi_{p,q}\colon\RR^{r\times s}\to\RR^{r\times s}$ by
\begin{equation}\label{eq:correction-maps}
 \begin{aligned}
 p_i(\xi)&=e^{\mathbf{i} h_{p,q}(\xi)}p_i e^{-\mathbf{i} h_{p,q}(\xi)}
                 \qquad(1\leqslant i\leqslant r),\\
 \Phi_{p,q}(\xi)
 &=\bigl(\tau_\omega (p_i(\xi)q_j)-\tau_\omega (p_i)\tau_\omega (q_j)\bigr)_{
                 1\leqslant i\leqslant r,\,1\leqslant j\leqslant s}.
 \end{aligned}
\end{equation}
Then 
 $\Phi_{p,q}(\xi)=0$ is equivalent to
$e^{\mathbf{i}h_{p,q}(\xi)}\cA_0e^{-\mathbf{i}h_{p,q}(\xi)}
\perp_{\tau_\omega}\cB_0$.
Since $h_{p,q}(\xi)$ is self-adjoint, each $p_i(\xi)$ is a projection with
$\tau_\omega(p_i(\xi))=\tau_\omega(p_i)$.  
Define the Gram operator
$G_{p,q}\colon\RR^{r\times s}\to\RR^{r\times s}$ for $p$ and $q$ by
\begin{equation}\label{eq:gram-operator}
 \langle G_{p,q}\xi,\eta\rangle
 =\tau_\omega \bigl(h_{p,q}(\xi)h_{p,q}(\eta)\bigr)
 \qquad (\xi,\eta\in\RR^{r\times s}).
\end{equation}
Equipping $\RR^{r \times s}$ with its Euclidean inner
product and norm $|\xi|$, we have
\begin{equation*}
 \langle G_{p,q}\xi,\xi\rangle=\|h_{p,q}(\xi)\|_2^2\geqslant0,
 \qquad
 \|h_{p,q}(\xi)\|\leqslant2\sum_{i,j}|\xi_{ij}|
 \leqslant2\sqrt{rs}|\xi|.
\end{equation*}

By \Cref{lem:fd-lifts}, we can choose unital trace-preserving $^*$-embeddings
\mbox{$\varphi_n^0:\cA_0\to\cM$} and $\psi_n^0:\cB_0\to\cM$ such that
\begin{equation*}
  a=(\varphi_n^0(a))_\omega
  \quad(a\in\cA_0),\qquad
  b=(\psi_n^0(b))_\omega
  \quad(b\in\cB_0).
 \end{equation*}
Write $p_{i,n}=\varphi_n^0(p_i)$ and $q_{j,n}=\psi_n^0(q_j)$.
These are projections satisfying
\begin{equation*}
 p_i=(p_{i,n})_\omega,\qquad q_j=(q_{j,n})_\omega,\qquad
 \tau(p_{i,n})=\tau_{\omega}(p_i),\qquad\tau(q_{j,n})=\tau_{\omega}(q_j)
 \quad(n\geqslant 1).
\end{equation*}
For each $n \geqslant 1$, define
\begin{equation*}
  \begin{aligned}
  h_{p,q}^{(n)}(\xi)
  &=\mathbf{i}\sum_{i=1}^r\sum_{j=1}^s
    \xi_{ij}[p_{i,n},q_{j,n}],
  &
  p_{i,n}(\xi)
  &=e^{\mathbf{i}h_{p,q}^{(n)}(\xi)}p_{i,n}
    e^{-\mathbf{i}h_{p,q}^{(n)}(\xi)},\\
  \Phi_{p,q}^{(n)}(\xi)
  &=\bigl(\tau(p_{i,n}(\xi)q_{j,n})
    -\tau_\omega(p_i)\tau_\omega(q_j)\bigr)_{i,j},
  &
  \langle G_{p,q}^{(n)}\xi,\eta\rangle
  &=\tau\bigl(h_{p,q}^{(n)}(\xi)h_{p,q}^{(n)}(\eta)\bigr).
  \end{aligned}
  \end{equation*}
Under this definition, $\Phi_{p,q}^{(n)} (\xi)=0$ is equivalent to orthogonality of the images of
$\operatorname{Ad}(e^{\mathbf{i}h_{p,q}^{(n)}(\xi)})\circ\varphi_n^0$
and $\psi_n^0$. 
 Denote
by $D\Phi_{p,q}^{(n)}(\xi)$   the derivative of $\Phi_{p,q}^{(n)}$ at $\xi$.
The following lemma provides a uniform derivative estimate
for the finite-dimensional lifting argument.

\begin{lemma}\label{lem:gram-derivative}
  With the notation above, we have $D\Phi_{p,q}^{(n)}(0)=G_{p,q}^{(n)}$ for each $n \geqslant 1$ and
  \begin{equation}\label{eq:derivative-bound}
   \|D\Phi_{p,q}^{(n)}(\xi)-G_{p,q}^{(n)}\|\leqslant C|\xi|
   \qquad(\xi\in\RR^{r\times s}),
  \end{equation}
  where $C=16(rs)^{\frac{3}{2}}$.
  \end{lemma}
  
  \begin{proof}
  For each $\xi\in\RR^{r\times s}$,
  define bounded linear maps $\delta_n(\xi)$, $\cE_n(\xi)$ on $\cM$ by
  \begin{equation*}
   \delta_n(\xi)x=\mathbf{i}[h_{p,q}^{(n)}(\xi),x],\qquad
   \cE_n(\xi)=e^{\delta_n(\xi)}.
  \end{equation*}
  Then
  \begin{equation*}
   \cE_n(\xi)x=e^{\mathbf{i}h_{p,q}^{(n)}(\xi)}x e^{-\mathbf{i}h_{p,q}^{(n)}(\xi)},
   \qquad
   \|\delta_n(\xi)\|\leqslant2\|h_{p,q}^{(n)}(\xi)\|
   \leqslant4\sqrt{rs}|\xi|.
  \end{equation*}
  The exponential derivative formula
  \begin{equation*}
   \left.\frac{d}{d\theta}e^{A+\theta B}\right|_{\theta=0}
   =\int_0^1e^{(1-t)A}Be^{tA}\,dt
  \end{equation*}
  for bounded operators $A,B$, applied to
  $A=\delta_n(\xi)$ and $B=\delta_n(\eta)$, gives
  \begin{equation*}
   D\cE_n(\xi)\eta
   =\left.\frac{d}{d\theta}\cE_n(\xi+\theta\eta)\right|_{\theta=0}
   =\int_0^1
   \cE_n((1-t)\xi)\delta_n(\eta)\cE_n(t\xi)\,dt.
  \end{equation*}
  At $\xi=0$, we have $\cE_n(0)=I$ and
  $D\cE_n(0)\eta=\delta_n(\eta)$.
  Since $p_{i,n}(\xi)=\cE_n(\xi)p_{i,n}$, we have
  \begin{equation*}
   \langle D\Phi_{p,q}^{(n)}(0)\eta,\zeta\rangle
   =\tau\bigl(h_{p,q}^{(n)}(\eta)h_{p,q}^{(n)}(\zeta)\bigr)
   =\langle G_{p,q}^{(n)}\eta,\zeta\rangle.
  \end{equation*}
  Hence $D\Phi_{p,q}^{(n)}(0)=G_{p,q}^{(n)}$.
  
  For fixed $\xi$, the linearity of $\delta_n$ implies
  $\cE_n(\theta\xi)=e^{\theta\delta_n(\xi)}$.
  Differentiating yields
  \begin{equation*}
   \frac{d}{d\theta}\cE_n(\theta\xi)
   =\cE_n(\theta\xi)\delta_n(\xi),
  \end{equation*}
  and hence
  \begin{equation*}
   \cE_n(t\xi)-I
   =\int_0^t\cE_n(\theta\xi)\delta_n(\xi)\,d\theta.
  \end{equation*}
  Since $\cE_n(\theta\xi)$ is an isometry, we have
  \begin{equation*}
   \|\cE_n(t\xi)-I\|\leqslant t\|\delta_n(\xi)\|
   \qquad(0\leqslant t\leqslant1).
  \end{equation*}
  Thus
  \begin{equation*}
   \begin{aligned}
   D\cE_n(\xi)\eta-\delta_n(\eta)=\int_0^1
     (\cE_n((1-t)\xi)-I)\delta_n(\eta)\cE_n(t\xi)\,dt
   +\int_0^1\delta_n(\eta)(\cE_n(t\xi)-I)\,dt.
   \end{aligned}
  \end{equation*}
  Consequently,
  \begin{equation*}
   \begin{aligned}
   \|D\cE_n(\xi)\eta-\delta_n(\eta)\|
   \leqslant\int_0^1((1-t)+t)
         \|\delta_n(\xi)\|\,\|\delta_n(\eta)\|\,dt
   =\|\delta_n(\xi)\|\,\|\delta_n(\eta)\|.
   \end{aligned}
  \end{equation*}
Since  $\|p_{i,n}\| \leqslant 1$ and $\|q_{j,n}\|\leqslant1$, we obtain
  \begin{equation*}
   \begin{aligned}
   |(D\Phi_{p,q}^{(n)}(\xi)-G_{p,q}^{(n)})\eta|
   \leqslant\sqrt{rs}\,
         \|D\cE_n(\xi)\eta-\delta_n(\eta)\|\leqslant16(rs)^{\frac{3}{2}}|\xi|\,|\eta|.
   \end{aligned}
  \end{equation*}
  This proves \eqref{eq:derivative-bound}.
  \end{proof}

Now we show that positive definiteness of the Gram operator
ensures the existence of orthogonal trace-preserving lifts.
  
\begin{lemma}\label{lem:positive-gram}
  Let 
   $\cA_0$ and $\cB_0$ be nontrivial finite-dimensional von Neumann
  subalgebras of $\cM^\omega$ with
  $\cA_0\perp_{\tau_\omega}\cB_0$.
  Suppose that $p=(p_i)_{i=1}^r$ and $q=(q_j)_{j=1}^s$
are tuples of projections such that
$\{1,p_1,\ldots,p_r\}$ and $\{1,q_1,\ldots,q_s\}$
form linear bases of $\cA_0$ and $\cB_0$ respectively,
  and   $G_{p,q}$ is positive definite.
  Then there are unital trace-preserving $^*$-embeddings
  $\varphi_n\colon \cA_0\to\cM$ and $\psi_n\colon \cB_0\to\cM$
    satisfying
  \begin{equation*}
    a=(\varphi_n(a))_\omega \,\,(a\in\cA_0),
   \quad
   b=(\psi_n(b))_\omega \,\,(b\in\cB_0),
  \quad 
   \varphi_n(\cA_0)\perp_\tau\psi_n(\cB_0) \,\,(n\geqslant 1).
  \end{equation*}
  \end{lemma}
  
  \begin{proof} 
  Since 
  $\Phi_{p,q}^{(n)} (0)\xlongrightarrow{\omega} 0$ and $h_{p,q}(\xi)=(h_{p,q}^{(n)} (\xi))_\omega$ for each $\xi$,
  we have 
  \begin{equation*}
   \lim_{n\to\omega}\langle G_{p,q}^{(n)} \xi,\eta\rangle
   =\langle G_{p,q}  \xi,\eta\rangle
   \qquad(\xi,\eta\in\RR^{r\times s}).
  \end{equation*}
  Since $\RR^{r\times s}$ is finite-dimensional,
  $\|G_{p,q}^{(n)} -G_{p,q} \| \xlongrightarrow{\omega} 0$.
  Let $C$ be the constant in \Cref{lem:gram-derivative} and
  $\lambda=\lambda_{\min}(G_{p,q})>0$.
  On the set
  \begin{equation*}
   E_0=\{n\in\NN:\|G_{p,q}^{(n)} -G_{p,q} \|<\frac{\lambda}{2}\}\in\omega,
  \end{equation*}
  we have
  \begin{equation*}
   \langle G_{p,q}^{(n)} \xi,\xi\rangle\geqslant\frac{\lambda}{2}|\xi|^2,
   \qquad \| (G_{p,q}^{(n)})^{-1}\|\leqslant \frac{2}{\lambda}:=L.
  \end{equation*}
  
  Set $\rho_n=2L|\Phi_{p,q}^{(n)}(0)|$ for $n\in E_0$
and $\rho_n=0$ for $n\notin E_0$.
For $n\in E_0$, define
\[
T_n(\xi)=\xi-(G_{p,q}^{(n)})^{-1}\Phi_{p,q}^{(n)}(\xi).
\]
Since 
   $\rho_n\xlongrightarrow{\omega} 0$, we have 
  $$E:=\big\{n\in E_0:LC\rho_n\leqslant \frac{1}{2} \big\} \in \omega.$$
  For $n\in E$ and $|\xi|\leqslant\rho_n$,
  \begin{equation*}
   \|DT_n(\xi)\|
   =\|(G_{p,q}^{(n)})^{-1}(G_{p,q}^{(n)}-D\Phi_{p,q}^{(n)}(\xi))\|
   \leqslant LC|\xi|\leqslant\frac12.
  \end{equation*}
  For $\xi,\eta\in\overline B(0,\rho_n):=\{ \xi \in \RR^{r \times s}:|\xi| \leqslant \rho_n\}$,
  integrating along  the line segment
  \[
  \gamma(t)=\eta+t(\xi-\eta) \qquad (0\leqslant t\leqslant1),
  \]
  we obtain
  \begin{align*} 
  |T_n(\xi)-T_n(\eta)| &\leqslant \int_0^1 |DT_n(\gamma(t))\cdot  (\xi-\eta)| dt \\
  &\leqslant \int_0^1 \| DT_n(\gamma(t)) \| \cdot |\xi-\eta|dt \leqslant \frac{1}{2}|\xi-\eta|.
  \end{align*} 
  Hence  
  $T_n$ is a $\frac{1}{2}$-Lipschitz map 
  on  
  $\overline B(0,\rho_n)$
  when $n \in E$.
  Since $$|T_n(0)|\leqslant L|\Phi_{p,q}^{(n)}(0)|=\frac{\rho_n}{2},$$ we also have
  \begin{equation*}
   |T_n(\xi)|\leqslant|T_n(\xi)-T_n(0)|+|T_n(0)|
   \leqslant\tfrac12|\xi|+\tfrac12\rho_n\leqslant\rho_n.
  \end{equation*}
  By Banach's fixed-point theorem, $T_n$ has a unique fixed point
  $\xi_n\in\overline B(0,\rho_n)$, and hence $\Phi_{p,q}^{(n)}(\xi_n)=0$.

  For $n \in E$,
  define unital trace-preserving  $^*$-embeddings
  \begin{equation*}
   \varphi_n=\operatorname{Ad}(e^{\mathbf{i}h_{p,q}^{(n)}(\xi_n)})\circ\varphi_n^0,
   \qquad \psi_n=\psi_n^0.
  \end{equation*}
  It follows from  $\Phi_{p,q}^{(n)}(\xi_n)=0$ that  
  $\varphi_n(\cA_0) \perp_\tau \psi_n(\cB_0)$.
  Moreover,
  for each  $n\in E$ and 
   \(a\in\cA_0\),
  \[
  \begin{aligned}
  \|\varphi_n(a)-\varphi_n^0(a)\|_2
  &=\|e^{\mathbf{i}h_{p,q}^{(n)}(\xi_n)}\varphi_n^0(a)e^{-\mathbf{i}h_{p,q}^{(n)}(\xi_n)}-\varphi_n^0(a)\|_2 \leqslant 2\|a\|\,\|e^{\mathbf{i}h_{p,q}^{(n)}(\xi_n)}-1\|\\
  &\leqslant  2\|a\|\,\|h_{p,q}^{(n)}(\xi_n)\|\leqslant 4 \|a \| \sqrt{rs}\,\rho_n.
  \end{aligned}
  \]
  Also,
  $\psi_n(b)=\psi_n^0(b)$
  for each $n \in E$ and $b \in \cB_0$.
  
  For $n\notin E$, fix a   unital trace-preserving $^*$-embedding
  \[
   \theta\colon
   \bigl(\cA_0\otimes \cB_0,\,
         \tau_\omega|_{\cA_0}\otimes\tau_\omega|_{\cB_0}\bigr)
   \longrightarrow(\cM,\tau),
  \] and define unital  trace-preserving $^*$-embeddings
  \[
   \varphi_n(a)=\theta(a\otimes1),\qquad
   \psi_n(b)=\theta(1\otimes b)
   \qquad(a\in \cA_0,\ b\in \cB_0).
  \]
  Then 
  \[
   \tau\bigl(\varphi_n(a)\psi_n(b)\bigr)
   =\tau\bigl(\theta(a\otimes b)\bigr)
   =\tau_\omega(a)\tau_\omega(b)=\tau(\varphi_n(a))\tau(\psi_n(b)).
  \] 
  Since $E\in\omega$ and  $\rho_n\xlongrightarrow{\omega}0$, we have 
  \[
   a=(\varphi_n(a))_\omega,\qquad
   b=(\psi_n(b))_\omega
   \qquad(a\in \cA_0,\ b\in \cB_0).
  \]
  Hence $\varphi_n$ and $\psi_n$ have all the required properties.
  \end{proof}

  We write  
  $\mathring{x}=x-\tau(x)1$ for each $x \in \cM$, 
  and   $\mathring{\cP}=\cP\cap\ker\tau$
  for each unital $^*$-subalgebra $\cP\subseteq\cM$.
The next lemma allows us to make
  the Gram operator positive definite by an arbitrarily small
  unitary conjugation while preserving orthogonality.

\begin{lemma}\label{lem:free-gram}
 Let $(\cM,\tau)$ be a type $\mathrm{II}_1$ factor, 
  let
  $p=(p_i)_{i=1}^r$ and $q=(q_j)_{j=1}^s$ be tuples of
  projections in $\cM$ such that  $\{1,p_1,\ldots,p_r\}$ and
  $\{1,q_1,\ldots,q_s\}$ are linearly independent, and let 
$v\in\cM$ be a trace-zero symmetry such that $W^*(v)$ and \mbox{$\cP=C^*(p,q)$} are free.
For $t\in(0,\frac{\pi}{2})$, set $w_t=e^{\mathbf{i}tv}$,
$q_j^t=w_tq_jw_t^*$, and $q^t=(q_j^t)_{j=1}^s$.

\begin{enumerate}[label=(\arabic*), leftmargin=*]
  \item[$(1)$] Define $\Phi_{p,q}$ as in
  Equation~\eqref{eq:correction-maps}, with
  $(\cM^\omega,\tau_\omega)$ replaced by $(\cM,\tau)$.
  Then
  \[
    \Phi_{p,q^t}(0)=\cos^2 t\,\Phi_{p,q}(0).
  \]

  \item[$(2)$] $G_{p,q^t}$ is positive definite.
\end{enumerate}
\end{lemma}

\begin{proof}
\textup{(1)} 
Since $v$ is free independent to $\cP$, we get
$$\tau(avb)=\tau(abv)=0,\   \ \tau(avbv)=\tau(a)\tau(b), \    \  a,b\in\cP.$$ 
It is easy to check that  
$\Phi_{p,q^t}(0)=\cos^2t\,\Phi_{p,q}(0)$
for each $t\in(0,\frac{\pi}{2})$.

\textup{(2)} Let $\cH$ be the closed linear span of
$\{xvyv\colon x,y\in\mathring{\cP}\}$ in $L^2(\cM,\tau)$ and $e_{\cH}$ the orthogonal projection from $L^2(\cM,\tau)$ onto $\cH$.
It is routine to check that 
\begin{equation*}
 \langle x'vy'v,xvyv\rangle_2
 =\tau(x'x^*)\tau(y'y^*)
 \qquad(x,y,x',y'\in\mathring{\cP}).
\end{equation*}
Thus $x\otimes y\mapsto xvyv$ is a norm-preserving map from
$\mathring{\cP}\otimes\mathring{\cP}$ to $\cH$.
By the freeness of $W^*(v)$ and $\cP$, we can show that
\begin{equation*}
 \cH\perp\cP,\qquad
 \cH\perp\{avb:a,b\in\cP\},\qquad
 \cH\perp\{vxvy:x,y\in\mathring{\cP}\}.
\end{equation*}
Notice that
\begin{equation*}
 \begin{aligned}
 \mathbf{i}[p_i,q_j^t]
 =\mathbf{i}\cos^2t[\mathring p_i,\mathring q_j]
 +\mathbf{i}\sin^2t
   (\mathring p_iv\mathring q_jv-v\mathring q_jv\mathring p_i)
 -\sin t\cos t
   [\mathring p_i,v\mathring q_j-\mathring q_jv],
 \end{aligned}
\end{equation*}
hence
$e_{\cH}(\mathbf{i}[p_i,q_j^t])
=\mathbf{i}\sin^2t\,\mathring p_iv\mathring q_jv$.
Since both the families $\{\mathring p_i\}_{i=1}^r$ and
$\{\mathring q_j\}_{j=1}^s$ are
linearly
independent,
we can deduce that for each nonzero $\xi\in\RR^{r\times s}$,
\begin{equation*}
 \begin{aligned}
 \langle G_{p,q^t}\xi,\xi\rangle
 \geqslant\|e_{\cH}h_{p,q^t}(\xi)\|_2^2=\sin^4t
 \left\|\sum_{i=1}^r\sum_{j=1}^s
       \xi_{ij}\mathring p_i\otimes\mathring q_j\right\|_2^2>0.
 \end{aligned}
\end{equation*}
Hence $G_{p,q^t}$ is positive definite.
\end{proof}
 
Now we prove the finite-dimensional lifting result. 

\begin{proposition}\label{prop:fd}
Let $(\cM,\tau)$ be a 
type $\mathrm{II}_1$ factor.
Suppose that $\cA_0,\cB_0\subseteq\cM^\omega$ are finite-dimensional von Neumann
subalgebras and $\cA_0\perp_{\tau_\omega}\cB_0$. Then there are unital
trace-preserving $^*$-embeddings $\varphi_n:\cA_0\to\cM$ and
$\psi_n:\cB_0\to\cM$ 
satisfying
\begin{equation*}
  a=(\varphi_n(a))_\omega \,\,(a\in\cA_0),
 \quad
 b=(\psi_n(b))_\omega \,\,(b\in\cB_0),
\quad 
 \varphi_n(\cA_0)\perp_\tau\psi_n(\cB_0) \,\,(n\geqslant 1).
\end{equation*}
\end{proposition}

\begin{proof}
By \Cref{lem:fd-lifts}, it suffices to
prove the case when $\cA_0$ and $\cB_0$
are nontrivial.
Suppose that $p=(p_i)_{i=1}^r$ and $q=(q_j)_{j=1}^s$
are tuples of projections such that
$\{1,p_1,\ldots,p_r\}$ and $\{1,q_1,\ldots,q_s\}$
form linear bases of $\cA_0$ and $\cB_0$ respectively.
By \Cref{thm:popa},
we can  choose a trace-zero symmetry $v$ free independent to
$C^*(\cA_0,\cB_0)$. For $0<t<\frac{\pi}{2}$, set $w_t=e^{\mathbf{i}tv}$,
$q_j^t=w_tq_jw_t^*$,  $q^t=(q_j^t)_{j=1}^s$ and
$\cB_0^t=w_t\cB_0w_t^*$.
By \Cref{lem:free-gram}, we have
$\cA_0\perp_{\tau_\omega}\cB_0^t$ and $G_{p,q^{t}}$ 
is positive definite. 

Let $S\subseteq(\cM^{r+s})_1$
be the set  of all tuples
\begin{equation*}
 (\varphi(p_1),\ldots,\varphi(p_r),
  \psi(q_1),\ldots,\psi(q_s)),
\end{equation*}
where $\varphi:\cA_0\to\cM$ and $\psi:\cB_0\to\cM$
run over all unital
trace-preserving $^*$-embeddings with
$\varphi(\cA_0) \perp_\tau \psi(\cB_0)$.
This set is nonempty by a unital trace-preserving $^*$-embedding
from
$\cA_0\otimes \cB_0$ to $\cM$. 
Set $S_n=S$ for every $n \geqslant 1$.
Fix $0<t<\frac{\pi}{2}$. By \Cref{lem:positive-gram}, there exist
unital trace-preserving $^*$-embeddings
\[
 \varphi_n^t:\cA_0\to\cM,\qquad
 \widetilde\psi_n^t:\cB_0^t\to\cM
\]
such that
\[
 a=(\varphi_n^t(a))_\omega\,\,(a\in\cA_0),\quad
 b=(\widetilde\psi_n^t(b))_\omega\,\,(b\in\cB_0^t),\quad  \varphi_n^t(\cA_0)\perp_\tau
 \widetilde\psi_n^t(\cB_0^t)
 \,\,(n\geqslant 1).
\]
Define
\[
 \psi_n^t(b)=\widetilde\psi_n^t(w_tbw_t^*)
 \qquad(b\in\cB_0).
\]
Then $\psi_n^t:\cB_0\to\cM$ is a unital trace-preserving
$^*$-embedding with
$\psi_n^t(\cB_0)=\widetilde\psi_n^t(\cB_0^t)$.
Hence, by the definition of $S_n$,
\[
 \bigl(\varphi_n^t(p_1),\ldots,\varphi_n^t(p_r),
       \psi_n^t(q_1),\ldots,\psi_n^t(q_s)\bigr)
 \in S_n
 \qquad(n\geqslant 1).
\]
Since
$p_i=
 (\varphi_n^t(p_i))_\omega$,  
 $w_tq_jw_t^*=(\psi_n^t(q_j))_\omega$,
we obtain
\[
 (p_1,\ldots,p_r,w_tq_1w_t^*,\ldots,w_tq_sw_t^*)
 \in\prod_\omega S_n.
\]
Since  $\lim_{t \to 0}\| w_tq_jw_t^*-q_j \|_2=0$,
it follows from  \Cref{lem:closed-lifts} that 
\begin{equation*}
 (p_1,\ldots,p_r,q_1,\ldots,q_s)\in\prod_\omega S_n.
\end{equation*} 
As a result, there exist
unital trace-preserving $^*$-embeddings
$
 \varphi_n:\cA_0\to\cM$ and 
 $\psi_n:\cB_0\to\cM$
such that
\[
 \varphi_n(\cA_0)\perp_\tau\psi_n(\cB_0)
 \qquad(n\geqslant 1),
\]
and $p_i=(\varphi_n(p_i))_\omega$, 
$q_j=(\psi_n(q_j))_\omega$ for each $1 \leqslant i \leqslant r$ and 
$1 \leqslant j \leqslant s$.
Since $\{1,p_1,\ldots,p_r\}$ and $\{1,q_1,\ldots,q_s\}$
are linear bases,
$\varphi_n$ and $\psi_n$ satisfy all the required properties.
\end{proof}

Now we are in a position to prove  \Cref{thm:amenable} using
the finite-dimensional lifting result.

\begin{proof}[Proof of \Cref{thm:amenable}]
Since $\cA$ and $\cB$ are amenable and separable, 
there exist
 increasing sequences of
finite-dimensional subalgebras $\cA^{(k)}\subseteq \cA$ and
$\cB^{(k)}\subseteq \cB$ with weakly dense unions.
When $\cA$ (resp. $\cB$) is
abelian,   $\cA^{(k)}$ (resp. $\cB^{(k)}$) is also  abelian.
Since  $\cA^{(k)}\perp_{\tau_\omega} \cB^{(k)}$,
it follows from \Cref{prop:fd} that there exist
 trace-preserving $^*$-embeddings
\begin{equation*}
 \varphi_{k,n}\colon \cA^{(k)}\to \cM,\qquad
 \psi_{k,n}\colon \cB^{(k)}\to \cM
\end{equation*}
such that
\begin{align*}
  &a=(\varphi_{k,n}(a))_{n \to \omega}
  \quad(a\in\cA^{(k)}),\qquad \qquad
  b=(\psi_{k,n}(b))_{n \to\omega}
  \quad(b\in\cB^{(k)}),\\
  &\varphi_{k,n}(\cA^{(k)})
  \perp_\tau
  \psi_{k,n}(\cB^{(k)})\quad  (n,k \geqslant 1).
 \end{align*}

Since $\cA$ and $\cB$ are both separable, 
by Kaplansky's density theorem,  
we can choose sequences
\[
\{a_j\}_{j\geqslant1}\subseteq
\Bigl(\bigcup_k\cA^{(k)}\Bigr)_1,
\qquad
\{b_j\}_{j\geqslant1}\subseteq
\Bigl(\bigcup_k\cB^{(k)}\Bigr)_1
\]
that are  $2$-norm dense in $(\cA)_1$ and $(\cB)_1$
 respectively. 
 Passing to a common subsequence of the approximating algebras
and relabeling the corresponding $^*$-embeddings,
we may  assume that
$$a_1,\ldots,a_k\in \cA^{(k)},\quad b_1,\ldots,b_k\in \cB^{(k)} \qquad(k\geqslant1).$$
Write $a_j=(a_{j,n})_{n \to \omega}$ and
$b_j=(b_{j,n})_{n \to \omega}$,
with each $a_{j,n},b_{j,n} $ a contraction 
in $\cM$.
For $k\geqslant1$, define
\begin{equation*}
 F_k=\left\{n\in\NN\colon
 \max_{1\leqslant j\leqslant k}
 \|\varphi_{k,n}(a_j)-a_{j,n}\|_2<\frac1k,\quad 
 \max_{1\leqslant j\leqslant k}
 \|\psi_{k,n}(b_j)-b_{j,n}\|_2<\frac1k 
 \right\}.
\end{equation*}
Then $F_k\in\omega$. Set
$
 E_k=\{n:n\geqslant k \} \cap\bigcap_{\ell=1}^k F_\ell,
$
then $E_k \in \omega$.
Define $$k_n=\max(\{1\}\cup\{k: n\in E_k\}).$$
Then $k_n\xlongrightarrow{\omega}\infty$. Define
\begin{equation*}
 \cA_n=\varphi_{k_n,n}(\cA^{(k_n)}),\qquad
 \cB_n=\psi_{k_n,n}(\cB^{(k_n)}).
\end{equation*}
These are finite-dimensional von Neumann subalgebras with $\cA_n\perp_{\tau}\cB_n$.

Now we show that 
$a_j \in \prod_\omega \cA_n$,  $b_j\in\prod_{\omega}\cB_n$  for each  $j\geqslant1$.
Define $a'_{j,n}\in(\cA_n)_1$ by
\begin{equation*}
 a'_{j,n}=\begin{cases}
 \varphi_{k_n,n}(a_j),&k_n\geqslant j,\\
 0,&k_n<j.
 \end{cases}
\end{equation*}
Then $a'_{j,n} \in \cA_n$ for each $n\geqslant1$. When $n \in  E_j$,
 we have $k_n\geqslant j$ and $n\in E_{k_n} \subseteq F_{k_n}$, hence
\begin{equation*}
 \|a'_{j,n}-a_{j,n}\|_2< \frac{1}{k_n}\xlongrightarrow{\omega}0,
\end{equation*}
which implies that $a_j=(a'_{j,n})_{n \to \omega}\in\prod_{\omega}\cA_n$.
Similarly, 
 $b_j\in\prod_{\omega}\cB_n$ for each $j \geqslant 1$.

 Since $\prod_\omega\cA_n$ and $\prod_\omega\cB_n$ are von Neumann
subalgebras of $\cM^\omega$, their unit balls are closed in 
$2$-norm. Notice that 
\[
(\cA)_1
=\overline{\{a_j:j\geqslant1\}}^{\|\cdot\|_2}
\subseteq\Bigl(\prod_\omega\cA_n\Bigr)_1,
\qquad
(\cB)_1
=\overline{\{b_j:j\geqslant1\}}^{\|\cdot\|_2}
\subseteq\Bigl(\prod_\omega\cB_n\Bigr)_1.
\]
Hence
$
\cA\subseteq\prod_\omega\cA_n$,
$
\cB\subseteq\prod_\omega\cB_n$.

If $\cA$ (respectively, $\cB$) is abelian, then each
$\cA_n$ (respectively, $\cB_n$) is abelian by the
choices of
$\cA^{(k)},\cB^{(k)}$.
\end{proof}

The following corollary of \Cref{thm:amenable} answers
Problem~30 in \cite{KE25} affirmatively.

\begin{corollary}\label{cor:unitaries}
Suppose that $(\cM,\tau)$ is a type $\mathrm{II}_1$ factor,
 $u,v$ are independent unitary
operators in $\cM^{\omega}$.
Then there are finite-spectrum unitary
operators
$u_n, v_n \in \cM$  such that
\begin{equation*}
 u=(u_n)_{\omega},\qquad v=(v_n)_{\omega},\qquad W^*\left(u_n\right) \perp_\tau W^*\left(v_n\right) \quad (n \geqslant1).
\end{equation*}
\end{corollary}

\begin{proof}
Applying
\Cref{thm:amenable} to $\cA=W^*(u)$
and $\cB=W^*(v)$,
there exist  abelian finite-dimensional
algebras $\cA_n,\cB_n \subseteq \cM$, 
with  $\cA_n\perp_{\tau}\cB_n$,  such that
$\cA \subseteq \prod_{\omega} \cA_n$
and $\cB \subseteq \prod_{\omega} \cB_n$.
Write  $u=(a_n)_\omega$ with $a_n\in\cA_n$.
Since $\cA_n$ is finite-dimensional, the partial isometry
in the polar decomposition of $a_n$ extends to a unitary
$u_n\in\cA_n$, so that $a_n=u_n|a_n|$.
As $u$ is unitary, we have
\[
 \|u_n-a_n\|_2
 =\|1-|a_n|\|_2
 \leqslant\|1-a_n^*a_n\|_2
 \xlongrightarrow{\omega}0,
\]
hence $u=(u_n)_\omega$.
Similarly, there exist unitary 
operators $v_n\in\cB_n$
such that $v=(v_n)_\omega$.
Since $\cA_n$ and $\cB_n$ are finite-dimensional,
$u_n$ and $v_n$ 
 have finite spectrum.
It follows from $\cA_n\perp_{\tau}\cB_n$ 
that 
$u_n$ and $v_n$ are independent for each $n \geqslant 1$.
\end{proof}

\begin{remark}
  The same argument as in the proof of \Cref{thm:amenable}
  shows that the conclusion remains valid for the  
  ultraproduct of a sequence $\{(\cM_n,\tau_n)\}_{n \geqslant 1}$
  of type $\mathrm{II}_1$ factors,  
  with $\cA_n,\cB_n\subseteq\cM_n$ and
  $\cA_n\perp_{\tau_n}\cB_n$ for each $n\geqslant 1$.
\end{remark}

\section{Counterexamples for nonamenable subalgebras}\label{sec:counterexamples}
In this section, we prove \Cref{thm:counterexample}, showing that the
amenability assumption in \Cref{thm:amenable} is essential.
We construct orthogonal subalgebras of $\cM^\omega$ that admit no
orthogonal  lifts.

\begin{proof}[Proof of \Cref{thm:counterexample}]
  Since every projection in $\cM^\omega$ is 
  unitarily conjugate 
to a projection in the diagonal copy of $\cM$,
and simultaneous
unitary conjugation preserves orthogonality and the existence
of orthogonal lifts, we may assume that $q\in\cM$.
Set $\cD=W^*(q)$.
We first construct a trace-preserving $^*$-embedding
$\alpha:\cN\to\cM^\omega$ such that $\alpha(\cN)$ and $\cD$
are orthogonal with respect to $\tau_\omega$ but admit
no orthogonal lifts.

Let $(\cR,\tau_{\cR})$ be the hyperfinite type $\mathrm{II}_1$
factor. By \cite[Theorem~4.7]{AKE21},
we can choose trace-preserving $^*$-embeddings
$\pi,\rho\colon \cN\to\cR^{\omega}$ such that
$\rho\ne(T_n)_{\omega}\circ\pi$ for every sequence of trace-preserving
unital completely positive maps $T_n\colon\cR\to\cR$.
Set  $p=1-q$, $s=\tau_\omega(p)$ and $t=\tau_\omega(q)$.
Choose trace-preserving $^*$-embeddings
\[
\iota_p:(\cR,\tau_{\cR})\longrightarrow \big( p\cM p,s^{-1}\tau  \big),
\qquad
\iota_q:(\cR,\tau_{\cR})\longrightarrow \big( q\cM q,{t}^{-1}\tau  \big),
\]
and define
\begin{equation*}
 \alpha=\iota_p^{\omega}\circ\pi+\iota_q^{\omega}\circ\rho,
 \qquad \cA_0=\alpha(\cN).
\end{equation*}
Then $\alpha$ is a trace-preserving $^*$-embedding and
\begin{equation*}
 \tau_\omega(q\alpha(x))=t\tau_\omega(\alpha(x))\qquad(x\in \cN).
\end{equation*}
Thus $\cA_0\perp_{\tau_\omega} \cD$.

Assume on the contrary that von Neumann subalgebras
\mbox{$\cA_n,\cD_n\subseteq\cM$} satisfy
$$\cA_0\subseteq\prod_\omega\cA_n,\qquad
\cD\subseteq\prod_\omega\cD_n,\qquad
\cA_n\perp_\tau\cD_n\quad(n\geqslant 1).$$
Choose projections $q_n\in\cD_n$ such that $q=(q_n)_\omega$.
Set
$p_n=1-q_n$, $s_n=\tau(p_n)$ and $t_n=\tau(q_n)$. 
Let $E_n$ denote the
trace-preserving conditional expectation 
from $\cM$ to $\cA_n$.
Then $E_n(p_n)=s_n1$ and $E_n(q_n)=t_n1$. 
Set $K=\iota_q^{-1}\circ E_{\iota_q(\cR)}\colon q\cM q\to\cR$,
$\lambda_n=\max\{s_n,s\frac{t_n}{t} \}$, and 
\[
 S_n(z)=\frac{1}{\lambda_n}K\bigl(qq_n E_n(p_n\iota_p(z)p_n)q_nq\bigr)
 \qquad(z\in\cR).
\]
Since
$p_n\iota_p(1)p_n\leqslant p_n$ and $s_n\leqslant\lambda_n$,
we can deduce that  $S_n$ is completely positive and $S_n(1)\leqslant1$. For $z\geqslant0$,
\begin{equation*}
 \begin{aligned}
 \tau_{\cR}(S_n(z))
 \leqslant\frac{\tau(q_n E_n(p_n\iota_p(z)p_n))}{t\lambda_n}
 =\frac{t_n\,\tau(p_n\iota_p(z)p_n)}{t\lambda_n}
 \leqslant\tau_{\cR}(z).
 \end{aligned}
\end{equation*}
By Kadison-Schwarz inequality \cite{Ka52}, $\| S_n(z) \|_2 \leqslant \| z \|_2$
for each $z \in \cR$. 

Let
$
E=(E_n)_\omega\colon
\cM^\omega\to \prod_\omega\cA_n
$
be the trace-preserving conditional expectation.
Since  $\cA_0\subseteq\prod_\omega\cA_n$ and $E(p)=s1$,
the bimodular property of $E$ yields
\[
E(p\alpha(x))=s\alpha(x)\qquad(x\in\cN).
\]
Using
\[
\iota_p^\omega(\pi(x))=p\alpha(x),\qquad
\iota_q^\omega(\rho(x))=q\alpha(x)q,
\]
and $\lambda_n\xlongrightarrow{\omega}s$, we obtain
  \begin{equation*}
    (S_n)_{\omega}(\pi(x))
    =\frac{1}{s}K^{\omega}(qE(p\alpha(x))q)=\rho(x)\qquad(x\in \cN).
   \end{equation*}
Let $d_n=1-S_n(1)$ and $\omega_n=\tau_{\cR}-\tau_{\cR}\circ S_n$.
Then $\|d_n\|_2 \xlongrightarrow{\omega}0$, and $\omega_n$ is positive with
$\omega_n(1)=\tau_{\cR}(d_n)$.
 Define
\begin{equation*}
 T_n(z)=
 \begin{cases}
 S_n(z)+\dfrac{\omega_n(z)}{\omega_n(1)}d_n,&d_n\ne0,\\[4pt]
 S_n(z),&d_n=0.
 \end{cases}
\end{equation*}
Then $T_n$ is a unital trace-preserving
completely positive  map on $\cR$.
For $\|z\|\leqslant1$, we have
\[
\|T_n(z)-S_n(z)\|_2\leqslant\|d_n\|_2
\xlongrightarrow{\omega}0,
\]
hence $\rho=(T_n)_{\omega}\circ\pi$, a contradiction.
Therefore $\alpha(\cN)$ and $\cD$
 admit
no orthogonal lifts.

By \cite[Corollary~0.2]{Po14}, there exists 
a unitary $u$ with $u\cD u^*=\cD$, 
such that  
$uW^*(\cA_0,\cD)u^*$  
and $\cB$ are free with amalgamation over $\cD$,
 thus
\[
E_{\mathcal D}(ab)
=E_{\mathcal D}(a)E_{\mathcal D}(b)
\qquad
(a\in u\mathcal A_0u^*,\, b\in\mathcal B).
\]
Let $\cA=u\cA_0u^*$. Since $E_\cD(a)=\tau_\omega(a)1$ for each  $a\in \cA$, we have
\begin{equation*}
 \tau_\omega(ab)=\tau_\omega(E_\cD(a)E_\cD(b))=\tau_\omega(a)\tau_\omega(b)
 \qquad(a\in \cA,\ b\in \cB),
\end{equation*}
which implies that $\cA \perp_{\tau_\omega} \cB$.
Since $u\cD{}u^*=\cD$,  $\cA$ and $\cD$ admit no orthogonal  lifts.
\end{proof}

 \textbf{AI Statement}: During the preparation of this manuscript, 
 the authors used 
GPT 6  
for 
searching necessary references, 
exploring  proof strategies, examining mathematical arguments, 
and editing the language. 
 The authors 
take full responsibility for the  manuscript.

\end{document}